\documentclass{amsart}

\parskip=\smallskipamount

\newtheorem{theorem}{Theorem}[section]
\newtheorem{lemma}[theorem]{Lemma}

\newtheorem{proposition}[theorem]{Proposition}

\theoremstyle{definition}

\newtheorem{remark}[theorem]{Remark}

\newcommand{\C}{{\mathbb C}}
\newcommand{\N}{{\mathbb N}}

\newcommand{\R}{{\mathbb R}}

\newcommand{\cC}{\mathcal{C}}
\newcommand{\cF}{\mathcal{F}}

\numberwithin{equation}{section}

\begin{document}
\title[A contractible Levi-flat hypersurface]
{A contractible Levi-flat hypersurface 
which is a determining set for pluriharmonic functions}
\author{Franc Forstneri\v c}
\address{Institute of Mathematics, Physics and Mechanics, 
University of Ljubljana, Jadranska 19, 1000 Ljubljana, Slovenia}
\email{franc.forstneric@fmf.uni-lj.si}
\thanks{Supported by grants P1-0291 and J1-6173, Republic of Slovenia.}

%
%    General info
%
\date{February 22, 2005} 

\begin{abstract} 
We find a real analytic Levi-flat hypersurface in $\C^2$
containing a bounded contractible domain which is a determining 
set for pluriharmonic functions.
\end{abstract}

\maketitle

\section{The main result}
A real hypersurface $M$ in an $n$-dimensional complex manifold 
is {\em Levi-flat} if it is foliated by complex manifolds of 
dimension $n-1$; this {\em Levi foliation} is as smooth as 
$M$ itself according to Barrett and Forn\ae ss \cite{BF}.
If $M$ is real analytic, it is locally near every point 
defined by a {\em pluriharmonic function} 
$v$: $dd^c v=2i\partial\overline \partial v =0$. 
One might expect that an oriented, real analytic, Levi flat hypersurface 
admits a pluriharmonic defining function on any topologically simple 
relatively compact domain, perhaps under an additional 
analytic asumption such as the existence of a 
fundamental system of Stein neighborhoods (see e.g.\  Theorem 2 in
\cite{T}, p.\ 298). Here we show that, on the contrary, 
even a most simple domain in a real analytic Levi flat 
hypersurface may be a {\em determining set for pluriharmonic functions}.

\begin{theorem} 
\label{Maintheorem}
There exist an ellipsoid $B\subset \C^2$ and a 
real analytic, Levi-flat hypersurface $M\subset\C^2$ 
intersecting the boundary $bB$ transversely such that 
the Levi foliation of $M$ has trivial holonomy 
and $A=M\cap B$ satisfies the following:
\begin{itemize}
\item[(i)] $\overline {A}$ is diffeomorphic to the three-ball and 
admits a Stein neighborhood basis.  
\item[(ii)] Any real analytic function on $A$ which is 
constant on Levi leaves is constant. 
\item[(iii)] Any pluriharmonic function in a connected open
neighborhood of $A$ in $\C^2$ which vanishes on $A$ is identically zero.
\end{itemize}
\end{theorem}

The Levi foliation of $M$ in our proof is a {\em simple foliation} 
(\cite{God1}, p.\ 79) whose leaves are complex discs. 
Likely one can also obtain a similar example in the ball of $\C^2$. 
On the other hand, for any compact subset $A$ in a 
real analytic, simply connected Levi-flat hypersurface $M$ 
there is a {\em smooth} defining function $v$ for $M$ 
whose pluricomplex Laplacian $dd^c v$ is {\em flat} on $A$;
this suffices for the construction of Stein neighborhood basis of 
certain compact subsets of $M$ \cite {FL}.

We mention that D.\ Barrett gave an example of a 
{\em compact} real analytic Levi-flat hypersurface 
with trivial holonomy and without a global pluriharmonic 
defining function (Theorem 3 in \cite{Ba}). 
His example is the quotient of $S^1\times\C^*$
by $(\theta,z)\to (\phi(\theta),2z)$ where $\phi$  is a real analytic 
diffeomorphism of the circle $S^1$ which is topologically 
but not diffeomorphically conjugate to a rotation.

\section{A real analytic foliation of $\R^2$ without analytic 
first integrals}
Our construction of the hypersurface $M$ in theorem \ref{Maintheorem} 
is based on the following.

\begin{proposition}
\label{proposition}
Let $D$ be the open unit disc in $\R^2$.
There exists a real analytic foliation $\cF$ of\ $\R^2$ by closed 
lines such that any real analytic function on $D$ which is 
constant on every leaf of the restricted foliation $\cF|_D$ is constant. 
\end{proposition}

\begin{remark} 
While we cannot exclude the possibility 
that an example of this kind is contained in the vast 
literature on the subject, we could not find a precise 
reference in some of the standard sources concerning foliations
of the plane (\cite{HR}, \cite{Ha}, \cite{God1}, \cite{God2}, \cite{CC}). 
It is known that every smooth foliation  of $\R^2$ by lines has a 
global continuous first integral but in general not one of 
class $\cC^1$, not even in the analytic case (Wazewsky \cite{W});
however, there exists a smooth first integral without critical
points on any relatively compact subset (Kamke \cite{K}).
\end{remark}

\begin{proof} 
Let $(x,y)$ be  coordinates on $\R^2$. Define subsets
$E_1,E_2\subset \R^2$ by
$$
	E_1 = \{(x,y)\in \R^2\colon x<-1 \ {\rm or}\ y>0\}, \quad
	E_2 = \{(x,y)\in \R^2\colon x>1 \  {\rm or}\ y>0\}.
$$
Let $\cF_j$ denote the restriction of the foliation $\{y=c\}_{c\in\R}$
to $E_j$ $(j=1,2)$. Let $\psi$ be a real analytic orientation 
preserving diffeomorphism of the half line $(0,+\infty)$,
so $\lim_{t \downarrow 0}\psi(t)=0$. (We do not require
that $\psi$ extends analytically to a neighborhood of $0$.) Then 
$\phi(x,y)=(x,\psi(y))$ is a real analytic diffeomorphism
of the upper half plane $E_{1,2}= E_1\cap E_2=\{(x,y)\in \R^2\colon y>0\}$ 
onto itself which maps every leaf of $\cF_1|_{E_{1,2}}$ to 
a leaf of $\cF_2|_{E_{1,2}}$. Let $E$ be the quotient of the
topological (disjoint) sum $E_1\sqcup E_2$ obtained by identifying a point 
$(x,y)\in E_1$ belonging to $E_{1,2}$ with the point 
$\phi(x,y)\in E_2$. The foliations $\cF_j$ $(j=1,2)$ amalgamate 
into a real analytic foliation $\cF$ on $E$. 

By construction
$E$ is a real analytic manifold homeomorphic to $\R^2$,
and hence there exists a real analytic diffeomorphism
of $E$ onto $\R^2$. (This follows in particular from 
the classification theorem for simply connected Riemann
surfaces.) We identify $E$ with $\R^2$ and denote the resulting
real analytic foliation of $\R^2$ by $\cF=\cF_\psi$. 
Let $\pi\colon \R^2 \to Q=\R^2/\cF$ 
denote the projection onto the space of leaves. 
$Q$ admits the structure of a non-Hausdorff real analytic manifold 
such that $\pi$ is a real analytic submersion. (The real analytic
structure on $Q$ is obtained by declaring the restriction 
of $\pi$ to any local analytic transversal $\ell$ to $\cF$
to be a diffeomorphism of $\ell$ onto the open set $\pi(\ell)\subset Q$.
For the details see \cite{Ha}, \cite{HR}.) 
In our case $Q$ is the quotient of the topological 
sum $\R_1 \sqcup \R_2$ of two copies of 
the real axis obtained by identifying a point $t>0$ in $\R_1$ 
with the point $\psi(t)\in \R_2$ (no identifications are made 
for points $t\le 0$). The only pair of branch points 
in $Q$ (i.e., points without a pair of disjoint neighborhoods) 
are those corresponding to $0\in \R_1$ and $0\in \R_2$.

\begin{lemma}
\label{lemma1}
If $\psi$ is flat at origin 
$(\lim_{t\downarrow 0} \psi^{(k)}(t)=0$ for $k\in \N)$
then every real analytic function on $\R^2$ which is 
constant on every leaf of $\cF_\psi$ is constant. 
\end{lemma}

\begin{proof} 
A real analytic function $f$ on $\R^2$ which is constant on the 
leaves of the foliation $\cF_\psi$ is of the form $f=h\circ \pi$ for some 
real analytic function $h\colon Q\to \R$, where $Q$ is the space of leaves. 
From our construction of the foliation it follows that $h$ is given by
a pair of real analytic functions $h_j \colon \R\to \R$
$(j=1,2)$ satisfying $h_1(t)=h_2(\psi(t))$ for $t>0$. 
As $t\downarrow 0$, the flatness of $\psi$ at $0$ implies 
that the derivative $h_1'$ is flat at $0$. Hence $h_1$, and 
therefore also $h_2$, are constant.  
\end{proof}

Fix $\psi$ and consider the following pair of subsets 
of  $E_1$ resp.\ $E_2$: 
\begin{eqnarray*}
	D_1 &=& \{(x,y)\in \R^2 \colon -3<x<-2,\ -1< y <+2\}, \\
	D_2 &=& \{(x,y)\in\R^2  \colon 2<x<3,\  -1<y< \psi(2) \} \cup\\ 
	    &&  \cup\, \{(x,y)\in\R^2 \colon  -3<x<3,\ \psi(1)<y<\psi(2)\}. 
\end{eqnarray*}
Let $D$ be the quotient of the disjoint sum $D_1 \sqcup D_2$ 
obtained by identifying any point $(x,y) \in D_1$ such that $1<y<2$ 
with the point $\phi(x,y)=(x,\psi(y))\in D_2$. Clearly $D$ is a 
simply connected domain with compact closure in $E \simeq \R^2$,
and the space of leaves $Q_D = D/\cF$ is a  non-Haudorff manifold 
with a simple branch at $t=1\in\R_1$ resp.\ $\psi(1)\in \R_2$.

\begin{lemma}
\label{lemma2}
If $\psi$ is flat at the origin then every real analytic function 
$f$ on $D$ which is constant on every leaf of $\cF_\psi|_D$ is constant.
\end{lemma}

\begin{proof}
As in lemma \ref{lemma1} such $f$ is of the form
$f=h\circ \pi$ for some real analytic function
$h$ on $Q_D=D/\cF_\psi$.
Such $h$ is given by a pair of real analytic functions 
$h_1\colon (-1,2)\to\R$, $h_2\colon (-1,\psi(2))\to\R$
satisfying $h_1(t)=h_2(\psi(t))$ for $1< t< 2$.
By analyticity this  relation persists on the largest
interval on which both sides are defined, which is $(0,2)$.
By flatness of $\psi$ at $0$ we conclude as in lemma \ref{lemma1} that
$h_1$ and $h_2$ must be constant. 
\end{proof}

Let $\cF=\cF_\psi$ be the foliation of $\R^2$ constructed 
above with the diffeomorphism $\psi(t)=t e^{-1/t}$ of $(0,+\infty)$ 
(which is flat at $0$). Let $D\subset\subset \R^2$ satisfy the 
conclusion of lemma \ref{lemma2}. Choose a disc containing $D$;
clearly lemma \ref{lemma2} still holds for this disc, and by 
an affine change of coordinates on $\R^2$ we may assume this 
to be the unit disc. This completes the proof 
of proposition \ref{proposition}.
\end{proof}

\begin{remark}
Proposition \ref{proposition} holds for any foliation 
$\cF_\psi$ constructed above for which the diffeomorphism
$\psi$ of $(0,+\infty)$ is such that $h\circ \psi$ does not 
extend as a real analytic function to a neighborhood of $0$ for 
any real analytic function $h$ near $0$. An example is $t^\alpha$ for 
an irrational $\alpha>0$. 
The foliation of $\R^2$ determined by the algebraic 1-form 
$\omega= (\alpha-x)(1+x)dy - xdx$ has the space of leaves 
$\cC^1$-diffeomorphic to the `simple branch' $Q$ determined 
by $\psi(t)=t^\alpha$ (\cite{God2}, p.\ 120); hence it might be possible 
to find a disc $D\subset \R^2$ satisfying the proposition 
\ref{proposition} for this foliation. These examples indicate that 
a real analytic foliation of $\R^2$ only rarely admits 
real analytic first integrals on large compact subsets. 
\end{remark}

\section{Proof of theorem 1.1}
Let $\cF$ be a real analytic foliation of $\R^2$ furnished
by the proposition \ref{proposition}
such that any real analytic function on $D= \{x_1^2+x_2^2<1\}\subset \R^2$
which is constant on the leaves of $\cF|_D$ is constant.
Denote by $(x_1+iy_1,x_2+iy_2)$ the coordinates on $\C^2$
and identify $\R^2$ with the plane $\{y_1=0, y_2=0\} \subset \C^2$.
Complexifying the leaves of $\cF$ we obtain the Levi foliation  
of a closed, real analytic, Levi-flat hypersurface $M$ 
in an open tubular neighborhood $\Omega\subset \C^2$ of $\R^2$.
Set $B =\{x_1^2+x_2^2+ c(y_1^2+y_2^2)<1\}$ where $c>0$
is chosen sufficiently large such that $\overline B\subset \Omega$.
Note that $B\cap \R^2=D$. A generic choice of $c$ 
insures that $M$ intersects $bB$ transversely 
(since transversality holds along $bD \cap M$).
Set $A=M\cap B\subset\subset M$. If $B$ is sufficiently thin 
(which is the case if $c$ is sufficiently large) 
then clearly $\overline A$ is diffeomorphic to the closed ball
in $\R^3$. If a real analytic function $u\in \cC^\omega(A)$ 
is constant on every Levi leaf of $A$ then $u|_D$ is constant on 
every leaf of $\cF|_D$ and hence is constant. Thus $A$ satisfies 
property (ii) in theorem \ref{Maintheorem}.

The foliation $\cF$ of $\R^2$ is transversely orientable 
and hence admits a transverse real analytic vector field $\nu$. 
Its complexification is a holomorphic vector field $w$ 
in a neighborhood of $\R^2$ in $\C^2$ such that
$iw$ is transverse to $M$ in a neighborhood of 
$\overline B$, provided that $B$ is chosen sufficiently thin. 
Moving $M$ off itself to either side by a short time 
flow of $iw$ in a neighborhood of $\overline B$ 
we obtain thin neighborhoods of 
$\overline A$ with two Levi-flat boundary components; 
intersecting these with $rB$ for $r>1$ close to $1$ 
gives a fundamental system of Stein neighborhoods of $\overline A$. 

Suppose that $v$ is a real pluriharmonic function 
in a connected open neighborhood of $A$ such that $v|_A=0$. 
For every point $x\in A$ there is an open connected 
neighborhood $U_x\subset B$ and a pluriharmonic function $u_x$ on $U_x$, 
determined up to a real constant, such that
$u_x+iv$ is holomorphic on $U_x$. 
Since $A$ is contractible, $H^1(A,\R)=0$ and hence the 
collection $\{u_x\}_{x\in A}$ can be assembled into a pluriharmonic
function $u$ in a neighborhood of $A$ such that $u+iv$ is holomorphic. 
Since $v|_A=0$, $u$ is constant on every Levi leaf on $A$ 
and hence constant by property (ii) of $A$. Thus $v$ is constant 
and hence identically zero. This proves theorem \ref{Maintheorem}.

\smallskip
\it Acknowledgements. \rm
I wish to thank C.\ Laurent--Thi\'ebaut for having bro\-ught 
the question to my attention and  Nils \O vrelid for suggesting 
to consider complexified foliations.

\bibliographystyle{amsplain}

\end{document}